\documentclass{article}
 \usepackage{amsmath}
 \usepackage{amsfonts}
 \usepackage{amssymb}
\usepackage{latexsym}
\newtheorem{theorem}{Theorem}

\newtheorem{lemma}[theorem]{Lemma}

\newtheorem{remark}[theorem]{Remark}

\newtheorem{proposition}[theorem]{Proposition}
 \allowdisplaybreaks

\def\thebibliography#1{\section*{\centerline{References}\markboth
 {References}{References}}\list
{[\arabic{enumi}]}{\settowidth\labelwidth{[#1]}\leftmargin\labelwidth
 \advance\leftmargin\labelsep
 \usecounter{enumi}}
 \def\newblock{\hskip .11em plus .33em minus .07em}
 \setlength{\itemsep}{0pt plus 2pt \relax}
 \sloppy
 \sfcode`\.=1000\relax}

\begin{document}
\title{\bf  Sobolev Trace Inequalities}
\author{Young Ja Park}
\date{} 
\maketitle

\begin{abstract}
The existence of extremal functions for 
the Sobolev trace inequalities is studied using the concentration 
compactness theorem. The conjectured extremal, the function of 
conformal factor, is considered 
and is proved to be an actual extremal function  
with extra symmetry condition on functions. 
One of the limiting cases of the Sobolev trace inequalities 
is investigated and the best constant for this case is computed. 
\end{abstract}

\section{Introduction}      
\indent

The classical Sobolev inequalities on ${\bf R}^{n}$ and 
the Sobolev trace inequalities on ${\bf R}^{n+1}_{+}$ are given by 
\begin{eqnarray}
 \left(
    \int_{{\bf R}^{n}} |f(x)|^{s} dx
  \right)^{r/s} 
  \leq 
  c_{r,s} 
  \left(
    \int_{{\bf R}^{n}} 
         |\nabla f(x)|^{r} dx
  \right), \;\;\;\; 
\frac{1}{s} = \frac{1}{r} - \frac{1}{n}, \label{SOBO}
\end{eqnarray}
where $c_{r,s}$ is a positive constant independent of the 
function $f$, and  
\[
 \left(
    \int_{{\bf R}^{n}} |f(x)|^{q} dx
  \right)^{p/q} 
  \leq 
  A_{p,q} 
  \left(
    \int_{{\bf R}^{n+1}_{+}} 
         |\nabla u(x,y)|^{p} dxdy
  \right), \;\;\;
\frac{1}{q} = \frac{n+1}{np} - \frac{1}{n},
\]
where $u$ is an extension of $f$ to the upper half-space, and 
$A_{p,q}$ is a positive constant independent of the function $u$.  
In general, Sobolev inequalities provide estimates of lower order 
derivatives of a function in terms of its higher order derivatives. 
Recently, the importance of having the sharp form of the 
inequalities has been recognized. 
For example, the solution to the Yamabe problem turns out to 
depend on knowledge of the best constant of (\ref{SOBO}). 
In order to obtain the sharp form of inequalities, we often  
consider the variational problem associated with it. 
Then, we ask if an extremal function (a minimizer or maximizer) 
exists subject to some constraints. 
In fact, the question of existence of an extremal function of 
the inequality is directly related to that  of existence of 
a solution to the partial differential equation (Euler-Lagrange equation) 
corresponding to the variational problem. 

The sharp form of the Sobolev trace inequality 
for the case $p = 2$ and $n > 1$ is  
\[
 \left(
  \int_{{\bf R}^{n}} \!\!\! |f(x)|^{2n/(n-1)} dx 
  \right)^{\frac{(n-1)}{n}} 
 \!\!\!\!\! \leq \! 
  \frac{1}{\sqrt{\pi}}
  \frac{1}{n-1}
   \!  \left[ 
       \frac{\Gamma(n)}{\Gamma(n/2)}
     \right]^{\frac{1}{n}} 
 \!\! \left(
  \int_{{\bf R}^{n+1}_{+}} 
       \!\!  |\nabla u(x,y)|^{2} dxdy 
  \right),
\]
and extremal functions for this inequality are given by 
$f(x) = (1+|x|^{2})^{-(n-1)/2}$. 
Since this inequality is conformally invariant, the extremal 
function given above is unique up to a conformal automorphism. 
W. Beckner \cite{Beckner2} proved this by inverting 
the inequality to a fractional integral on the dual space 
and using a special case of the sharp Hardy-Littlewood-Sobolev inequality. 
Independently, J. Escobar \cite{Escobar} proved this 
by exploiting the conformal invariance of this inequality 
and using characteristics of an Einstein metric. 
He defined a new metric conformal to the Euclidean 
metric on the ball, and proved that the metric is, 
in fact, an Einstein metric based on the information obtained 
from the Euler-Lagrange equation of the inequality, 
which implies that the new metric has zero curvature 
with constant mean curvature on the boundary.

An extremal function for the Sobolev trace inequality 
for $1 < p < n+1$ 
\begin{eqnarray}
 \left(
    \int_{{\bf R}^{n}} |f(x)|^{q} dx
  \right)^{\frac{p}{q}} 
  \leq 
  A_{p,q} 
  \left(
    \int_{{\bf R}^{n+1}_{+}} 
         |\nabla u(x,y)|^{p} dxdy
  \right), \;\;\;
\frac{1}{q} = \frac{n+1}{np} - \frac{1}{n} \label{Sobtra}
\end{eqnarray}
was conjectured as the function of the form 
$f(x) = (1+|x|^{2})^{-(n+1-p)/2(p-1)}$. 

First, we will be concerned with the existence of extremal 
functions for the Sobolev trace inequalities when $1 < p < n+1$. 
In the study of the existence of extremal functions, 
a compactness problem arises when we deal with inequalities 
defined on the spaces which are invariant under dilations 
and translations. 
In the case of the Sobolev trace inequalities, 
this question can be put in the following context.  
Let ${\bf T}$ be the trace operator mapping 
$W^{1,p}({\bf R}^{n})$ to $L^{q}({\bf R}^{n})$ 
where $\frac{1}{q} = \frac{n+1}{np} - \frac{1}{n}$. 
Then ${\bf T}$ is a bounded linear operator. 
Now we consider the smallest positive constant 
$A_{p,q}$  with which the inequality (\ref{Sobtra}) holds 
for all $u$ in $W^{1,p}({\bf R}^{n})$ and we ask if the best 
constant $A_{p,q}$ is attained for some function $u$. 
The question concerning the constant $A_{p,q}$ is equivalent 
to the following minimization problem: 
\[
\inf \left\{ \int_{{\bf R}^{n+1}_{+}} 
         |\nabla u(x,y)|^{p} dxdy
 : u \in W^{1,p}({\bf R}^{n}), 
  \int_{{\bf R}^{n}} |f(x)|^{q} dx = 1 \right\},
\]
where $u$ is an extension of $f$ to the upper half-space. 
It is evident that (\ref{Sobtra}) remains unchanged if we replace 
$u$ by $\sigma^{-n/q}u(\cdot/\sigma)$ for $\sigma > 0$. 
This implies possible defects of compactness on 
minimizing sequences of the problem  
in the sense that if $u$ is a minimizer, then 
$u_{\sigma} = \sigma^{-n/q}u(\cdot/\sigma)$ will be another  
minimizer for each $\sigma$, and if we let $\sigma \to 0$ or 
$\sigma \to \infty$, then $(u_{\sigma})$ converges weakly to $0$ 
(which is certainly not a minimizer) and $(|u_{\sigma}|^{q})$ 
either converges weakly to a Dirac delta function 
as $\sigma \to 0$, or spreads out as $\sigma \to \infty$. 
Using the concentration compactness principle of P. L. Lions \cite{Lions2}, 
it is proved that 
any minimizing sequence of the variational problem of 
the Sobolev trace inequality is relatively compact in 
$L^{q}({\bf R}^{n})$ up to translations and dilations, 
and there exists an extremal function. 

We will look at the conjectured extremal function 
\cite{Escobar} for 
the Sobolev trace inequality. We can prove 
that this function is an actual extremal if we 
assume extra symmetry for the functions considered. 
In particular, we will consider a space of functions of 
conformal factor $[(1+y)^{2}+|x|^{2}]$, 
where $(x,y) \in {\bf R}^{n+1}_{+}$. 
Then, by simple argument, we can easily show that it is indeed 
a minimizer for the Sobolev trace inequality restricted on 
the functions of conformal factor.  

We will treat the Sobolev trace inequality for the case 
with $p = 1$ separately. The existence of the extremal function 
for this case is not guaranteed by the argument 
used for $p$ with $1 < p < n+1$. 
This case can be thought of as one of the limit cases of the 
inequality and is very closely related to the isoperimetric 
inequality. We will show that 
the extremal function does not exist for this particular case.      
The sharp constant will be computed using a rearrangement 
technique on the functions on ${\bf R}^{n+1}_{+}$.  

\section{Concentration compactness lemmas                 
and the existence of an extremal function}                

\indent

The Sobolev trace theorem tells us that there is a bounded 
linear operator from $W^{1,p}({\bf R}^{n+1}_{+})$ to 
$L^{q}({\bf R}^{n})$ (which is called the trace operator) 
where $1/q = (n+1)/np - 1/n$. 
This means that there exists a positive constant $C_{0}$ for 
which the following inequality holds for any 
$u \in W^{1,p}({\bf R}^{n+1}_{+})$:
\[
  \left(
    \int_{{\bf R}^{n}} |u(x,0)|^{q} dx
  \right)^{1/q} 
  \leq 
  C_{0} 
  \left(
    \int_{{\bf R}^{n+1}_{+}} |\nabla u(x,y)|^{p} dxdy
  \right)^{1/p}.
\]
The question we want to ask is whether there exists 
an extremal function for which the best constant is attained. 
To that end, we look at the following minimization problem: 
\begin{eqnarray}
\inf \! \left\{
 {\bf J}(u) \equiv \! \int_{{\bf R}^{n+1}_{+}}\!\! 
    |\nabla u(x,y)|^{p} dxdy : \!\! 
            \int_{{\bf R}^{n}} \!\!\! |u(x,0)|^{q} dx = 1, 
                    \; u \in W^{1,p}({\bf R}^{n+1}_{+}) 
    \right\}. \label{fine}
\end{eqnarray}
If an extremal function for (\ref{fine}) exists, then it must 
satisfy the following Euler-Lagrange equation: 
for a positive constant $C$,
\begin{eqnarray}
       \left. \begin{array}{ll}
       \mbox{div}(|\nabla u|^{p-2}\nabla u) = 0, & \hspace{1in} \mbox{on}\;\; 
                                     {\bf R}^{n+1}_{+}        \\
                                    \\ 
         |\nabla u|^{p-2} \frac{\partial u}{\partial y} +
               C |u|^{q-2}u    = 0, & \hspace{1in} \mbox{on}\;\; 
                         \partial {\bf R}^{n+1}_{+}
        \end{array} \right\}. 
\end{eqnarray}
Consider a minimizing sequence $(u_{k})$ for $(\ref{fine})$. 
From the trace theorem, we know that the infimum 
is finite and we denote it by ${\bf I}$. So we have 
\[ 
{\bf I} = \lim_{k \rightarrow \infty}    {\bf J}(u_{k}),
\] 
with
\vspace{-.15in}
\[
    \; u_{k} \in  W^{1,p}({\bf R}^{n+1}_{+}) 
    \;\:\mbox{and} \;\; 
    \int_{{\bf R}^{n}} |u_{k}(x,0)|^{q} dx = 1 \:\; 
    \mbox{for each} \;\: k.
\] 
Now $(u_{k})$ is a bounded sequence in $ W^{1,p}({\bf R}^{n+1}_{+})$ 
and in $L^{q}({\bf R}^{n})$. 
We can find a subsequence (which we will also denote $(u_{k})$) and 
$ u \in W^{1,p}({\bf R}^{n+1}_{+})$ such that 
$(u_{k})$ converges weakly to $u$ in $W^{1,p}({\bf R}^{n+1}_{+})$ 
and $(u_{k}(x,0))$ converges weakly to $u(x,0)$ in $L^{q}({\bf R}^{n})$. 
Since the integrand of ${\bf J}(\cdot)$ is convex, 
${\bf J}(\cdot)$ is lower semicontinuous and we have
\vspace{-.1in} 
\[ 
  {\bf J}(u) 
    \leq 
  \liminf_{k \rightarrow \infty} 
  {\bf J}(u_{k}) = {\bf I}
\]
and 
\vspace{-.15in}
\[
 \| u \| _{L^{q}({\bf R}^{n})} 
    \leq 
  \liminf_{k \rightarrow \infty} 
  \| u_{k} \| _{L^{q}({\bf R}^{n})} 
    = 1. 
\] 
If $\| u \| _{L^{q}({\bf R}^{n})} = 1$, then $u$ is a minimizer. 
So the real question is whether or not 
$\| u \| _{L^{q}({\bf R}^{n})} = 1$.
Since ${\bf J}(\cdot)$ and the $L^{q}({\bf R}^{n})$ norm  
are invariant under the translations and under the scaling
\[
\{ v(\cdot) \mapsto \sigma^{-\frac{n}{q}} 
v(\frac{\cdot}{\sigma}) \}
\]
for any $\sigma > 0$, we may be so unfortunate as to 
choose a minimizing sequence which has possibilities 
of failures of the compactness. 
But a good news is that we can design translations 
and dilations to avoid the failure of compactness 
by the concentration compactness theorem. 
The proof of the existence of an extremal function for 
the Sobolev trace inequality was sketched by P. L. Lions 
in his paper \cite{Lions2}. 
His proof is based on the concentration compactness theorem. 
We start by stating the concentration compactness lemmas. 
Hereafter $B_{r}(x)$ represents the ball centered at $x$ 
with radius $r$ in ${\bf R}^{n}$, ${\bf R}^{n+1}$, 
or ${\bf R}^{N}$, which will be clear in the context.  

\begin{lemma}[\hspace{-.02in}Concentration\hspace{-.02in}     
Compactness \hspace{-.035in}I\hspace{-.015in}]                
 \label{con-1}                                                
Let $(\rho_{k}\!)$ be a sequence in $L^{1}\!({\bf R}^{N}\!)$  
satisfying $\rho_{k} \geq 0$ in ${\bf R}^{N}$ and             
$\int_{{\bf R}^{N}} \rho_{k} dx = \lambda$                    
($\lambda$ fixed).                                            
Then there exists a subsequence $(\rho_{k_{j}})$ of           
$(\rho_{k})$ satisfying                                       
one of the following possibilities:\\                            
$(i)$ $(${\bf Compactness}$)$ there exists a sequence           
$(y_{j})$ in ${\bf R}^{N}$ so that for any $\varepsilon > 0 $ 
there exists $ R \in (0, \infty)$ such that                    
\[                                                     
   \int_{B_{R} ( y_{j})} \rho_{k_{j}}(x) dx \geq 
\lambda - \varepsilon. 
\]
In this case, 
 $(\rho_{k_{j}}(\:\cdot \: + y_{j}))$ is called 
${\bf tight}$.\\
$(ii)$ $(${\bf Vanishing}$)$ for any positive real number $R$, 
\[ 
   \lim_{j \to \infty} 
   \sup_{y \in {\bf R}^{N}} 
   \int_{B_{R}(y)} \rho_{k_{j}}(x) dx 
    = 0.  \]
$(iii)$ $(${\bf Dichotomy}$)$ there exists $\alpha \in (0, \lambda)$ 
such that for $\varepsilon > 0$, there exist $j_{0} \geq 1$ 
and sequences 
$(\eta_{j}),(\xi_{j}) \in L^{1}({\bf R}^{N})$ satisfying 
for $j \geq j_{0}$, 
\[
\| \rho_{k_{j}} - (\eta_{j} + \xi_{j}) \|_{L^{1}({\bf R}^{N})} < \varepsilon,
\]
\begin{eqnarray*}
&            & \left| \; \int_{{\bf R}^{N}} 
               \eta_{j}(x) dx - \alpha \; \right| 
               \leq \varepsilon,\;\;
    \left| \; \int_{{\bf R}^{N}} \xi_{j}(x) dx 
        - (\lambda - \alpha) \; 
    \right| 
        \leq \varepsilon, 
\end{eqnarray*}
and \hspace{.7in} $dist(supp\; \eta_{j}, 
                        supp\; \xi_{j}) 
               \to \infty $  
               \;\; as $j \to \infty$,  \\  
    where $dist(A,B) \equiv 
        \inf \{ d(a,b) \mid a \in A \;
    \mbox{and}\; b \in B \}$. 
\end{lemma}

\begin{lemma}[Concentration Compactness II]   \label{con-3}     
Let $\mu,\nu$ be two bounded nonnegative measures on           
${\bf R}^{N}$ satisfying for some constant $C > 0$             
\begin{eqnarray}                              
   \left(
      \int _{{\bf R}^{N}} |\varphi|^{q} d\nu
   \right)^{1/q} 
   \leq C 
   \left(
      \int _{{\bf R}^{N}} |\varphi|^{p} d\mu
   \right)^{1/p}, \;\;\; 
   \varphi \in C_0^{\infty} ({\bf R}^{N}) \label{concomII}
\end{eqnarray}
where $1 \leq p < q \leq \infty$. Then there exist 
an at most countable set ${\cal L}$, families  
$(x_{l})_{l \in {\cal L}}$ of distinct points in 
${\bf R}^{N}$, and $(\nu_{l})_{l \in {\cal L}}$ 
in $(0, \infty)$ such that 
\begin{eqnarray*}
   \nu = \sum_{l \in {\cal L}} \nu_{l} \delta_{x_{l}},\;\;\; 
   \mu \geq C^{-p} \sum_{l \in {\cal L}} \nu_{l}^{p/q} 
   \delta_{x_{l}}. 
\end{eqnarray*}
Thus, in particular, \hspace{.03in}     
$\sum _{l \in {\cal L}} \nu_{l}^{p/q} < \infty$. 
If, in addition, 
$\nu({\bf R}^{N})^{1/q} \geq 
C \mu({\bf R}^{N})^{1/p}$, 
then ${\cal L}$ reduces to a single point and 
$\nu = \gamma \; \delta_{x_{0}} = 
\gamma^{- p/q}C^{p} \mu$, for some 
$x_{0} \in {\bf R}^{N}$ and for some constant $\gamma \geq 0$. 
\end{lemma}      

\noindent

We want to prove that there exists a 
function for which the following infimum ${\bf I}$ is attained: 
\begin{eqnarray}
\inf  
    \left\{\!
              J(u) \equiv \!\! \int_{{\bf R}^{n+1}_{+}}  
            \!\! |\nabla u(x,y) |^{p} dxdy : \!\!\! 
           \int_{{\bf R}^{n}} \!\!\! |u(x,0)|^{q} dx = 1, 
                     u \in W^{1,p}({\bf R}^{n+1}_{+}) 
    \right\}.   \label{*1}
\end{eqnarray}
We will replace $W^{1,p}({\bf R}^{n+1}_{+})$ 
by $W^{1,p}({\bf R}^{n+1})$ without loss of generality. 
In this section, we will assume that $p > 1$ and this will ensure 
that $ p < q $, which we need to apply Lemma \ref{con-3}
in the proof of the following theorem. The case for $p = 1$ will 
be treated later separately. 

\begin{theorem} \label{2.2.1} 
Let $(u_{k})$ be a minimizing sequence of $(\ref{*1})$. 
Then there exist $(\sigma_{k})$ in $(0, \infty)$ 
and $ \; (w_{k})$ in ${\bf R}^{n}$ such that the new minimizing 
sequence $(\tilde{u}_{k})$ given by
\[
   \tilde{u}_{k}(x,y) \equiv
   \sigma_{k}^{- \frac{n}{q}} 
   u_{k} \left( \frac{x-w_{k}}{\sigma_{k}}, 
    \frac{y}{\sigma_{k}} \right),\quad 
   x \in {\bf R}^{n}, \;  
   y \in {\bf R}
\]
is relatively compact in 
$L^{q}({\bf R}^{n})$. In particular, 
$(\ref{*1})$ has a minimum. 
\end{theorem}    
\noindent {\bf Proof :}
Let $P_{k}(x,y) \equiv |\nabla u_{k}(x,y)|^{p} + 
                  |u_{k}(x,0)|^{q} \otimes \delta_{0}(y) 
                + |u_{k}(x,y)|^{\frac{n+1}{n}q}$. 
Then $P_{k} \geq 0$ and 
$\int_{{\bf R}^{n+1}} P_{k}(x,y) dxdy \rightarrow L 
\geq {\bf I} + 1$ by the Sobolev embedding theorem. 
The idea is to show that we can prevent vanishing and dichotomy 
occurring for this sequence of functions by judicial choice 
of dilations and translations, so that 
we conclude the claim of the theorem by Lemma \ref{con-1}. 
Consider the concentration function $Q_{k}$ of $P_{k}$ defined as 
\[ 
   Q_{k}(t) \equiv \sup_{(x,y) \in 
   {\bf R}^{n} \times {\bf R}} 
   \int_{B_{t}((x,y))} 
    P_{k}(w,s) dw ds \;\;\;\; 
    \mbox{for} \;\; t > 0.
\]    
Then $(Q_{k})$ is a sequence of non-decreasing 
continuous functions on ${\bf R}^{+}$. 
For $\sigma > 0$, consider the concentration 
function $Q_{k}^{\sigma}$ of 
\[
\tilde{P}_{k}(x,y) \equiv |\nabla \tilde{u}_{k}(x,y)|^{p} 
   + |\tilde{u}_{k}(x,0)|^{q} 
   \otimes \delta_{0}(y) + 
   |\tilde{u}_{k}(x,y)|^{\frac{n+1}{n}q},
\] 
where $\tilde{u}_{k}(x,y)$ is defined as in the 
statement of the theorem with $\sigma_{k} = \sigma$. 
Then we have $Q_{k}^{\sigma}(t) = Q_{k}(\frac{t}{\sigma})$. 
So, we see a chance of vanishing occur. 
In order to avoid that, we take a sequence 
$(\sigma_{k})$ of dilations so that 
$Q_{k}^{\sigma_{k}}(1) = \frac{1}{2}$. We can see that
\vspace{-.1in}
\[ 
   \lim_{k \rightarrow \infty} 
     \sup_{(x,y) \in {\bf R}^{n} \times {\bf R}} 
        \int _{B_{R}((x,y))} 
     \tilde{P}_{k}(w,s) dw ds \geq \frac{1}{2} \quad
     \;\: \mbox{for} \:\; R \geq 1 
\] 
since $Q_{k}^{\sigma_{k}}(t) \geq \frac{1}{2}$ for 
$t \geq 1$. We prevented vanishing occurring by the choice 
of dilations. We will denote the new minimizing sequence 
\[
\left(\sigma_{k}^{-\frac{n}{q}}u_{k}
(x/\sigma_{k}, y/\sigma_{k})  \right)
\]
by $(u_{k})$. Now we show that dichotomy does not occur. 

\begin{lemma}
The dichotomy does not occur.
\end{lemma}

\noindent {\bf Proof :}
Suppose it occurs. Then there exists 
$\lambda^{*} \in (0,L)$ such that for any 
$\varepsilon > 0$ there exist 
$(w_{k}, \tilde{w}_{k}) \in {\bf R}^{n} \times {\bf R}$ 
and $R_{k}$, $k = 0,1,2, \cdots$ with $R_{k} > R_0$ 
(for $k = 1,2, \cdots$) and $R_{k} \rightarrow \infty$ so that  
\begin{eqnarray*}
    \left|\;
     \lambda^{*} - \int_{B_{R_{0}}((w_{k}, 
        \tilde{w}_{k}))} P_{k}(x,y) dx dy \; \right| 
\!\!\! \!\! & < & \!\!\!\! \!\varepsilon ,  \hspace{1in}  \\
  \left|\;
     (L - \lambda^{*}) - \int_{[B_{R_{k}}((w_{k}, 
        \tilde{w}_{k}))]^{C}} P_{k}(x,y) dx dy \; \right| 
\!\!\! \!\! & < & \!\!\!\! \!\varepsilon ,  \hspace{1in}  \\
    \int_{R_{0} < |(x,y) - (w_{k}, 
        \tilde{w}_{k})| < R_{k}} P_{k}(x,y) dx dy 
 \!\!\!\!  & < &\!\!\!\!\varepsilon, \hspace{1in} \\
     \mbox{supp} \: [P_{k} \chi_{B_{R_{0}}((w_{k}, \tilde{w}_{k}))}] 
\!\!\!\!  & \subset &\!\!\!\!\! B_{R_{0}}((w_{k}, \tilde{w}_{k})),   \\
     \mbox{supp} \: [P_{k}(1 - \chi_{B_{R_{k}}((w_{k}, \tilde{w}_{k}))})] 
\!\!\!\!  & \subset &\!\!\!\!\! [B_{R_{k}}((w_{k}, \tilde{w}_{k}))]^C,   \\
   \mbox{dist}\!\left( \mbox{supp} 
    [P_{k} \chi_{B_{R_0}((w_{k}, \tilde{w}_{k}))}], 
                 \mbox{supp} 
    [P_{k}(1 - \chi_{B_{R_k}((w_{k}, \tilde{w}_{k}))})]\! 
          \right)   
 \!\!\!   \!  & & \!\!\!\! \\
 \geq \mbox{dist} \left( B_{R_0}((w_{k}, \tilde{w}_{k})), 
    [B_{R_k}((w_{k}, \tilde{w}_{k}))]^{C} \right) 
  \! \!\!\!   & \to &\!\!\!\! \infty \;\;\mbox{as}\;\; k \to \infty.
\end{eqnarray*}
Consider $\xi, \eta \in C_{b}^{\infty}({\bf R}^{n+1})$ 
satisfying $0 \leq \xi, \eta \leq 1$, and 
\begin{eqnarray*} 
   \xi(x,y) & = & \left\{ \begin{array}{ll} 
                         1  & \mbox{if $|(x,y)| \leq 1$} \\
             0  & \mbox{if $|(x,y)| \geq 2$,} 
                       \end{array} 
               \right. \\  
   \eta(x,y)& = & \left\{ \begin{array}{ll} 
                        0  & \mbox{if $|(x,y)| \leq \frac{1}{2}$} \\
            1  & \mbox{if $|(x,y)| \geq 1$. } 
                       \end{array}
               \right.  
\end{eqnarray*}
We may take  $R_1$ so that $4 R_{1} \leq R_{k}$ for $k = 2,3, \cdots$. 
Define  
\[
\xi_{k}(x,y) \equiv  \xi(\frac{x - w_{k}}{R_{1}}, 
        \frac{y - \tilde{w}_{k}}{R_{1}})\;\;\; \mbox{and} \;\;\;  
\eta_{k}(x,y) \equiv  \eta(\frac{x - w_{k}}{R_{k}}, 
        \frac{y - \tilde{w}_{k}}{R_{k}}).
\]
We look at the following quantity: for $k$ large enough,  
\begin{eqnarray*}   
{\bf M}& = & \!  \int_{{\bf R}^{n+1}} \!\!\!
                             |\nabla u_{k}|^{p}
                             dxdy - 
                          \int_{{\bf R}^{n+1}} \!\!\!
                             |\nabla (u_{k} \xi_{k})|^{p}
                             dxdy - 
                          \int_{{\bf R}^{n+1}} \!\!\!
                             |\nabla (u_{k} \eta_{k})|^{p}
                             dxdy \\
              & = & \!  \int_{B_{R_{k}} - B_{R_{1}}} \!\!\!\!\!\!\!\!
                             |\nabla u_{k}|^{p}dxdy - 
                          \int_{B_{2R_{1}} - B_{R_{1}}}\!\!\!\!\! \!\!\!\!\!\!
                             |\nabla (u_{k} \xi_{k})|^{p}dxdy - 
                          \int_{B_{R_{k}} - 
                                B_{\frac{1}{2}R_{k}}}\!\!\!\!\! \!\!\!\!\!\! 
                             |\nabla (u_{k} \eta_{k})|^{p}dxdy \\ 
              & \equiv & \!  {\bf M}_{1} - {\bf M}_{2} - 
                      {\bf M}_{3}.
\end{eqnarray*} 
First, we have 
\vspace{-.15in}
\[
{\bf M}_{1} \leq  \int_{R_{0} < 
                             |(x,y) - (w_{k}, \tilde{w}_{k})| 
                                  < R_{k}} P_{k}(x,y) 
                             dxdy < \varepsilon.
\]
Using H\"older's inequality and the Sobolev embedding theorem 
together with the assumptions in the beginning of the lemma, we show  
\begin{eqnarray*}
{\bf M}_{2}^{1/p}  
              & \leq & 
                          \left(\int_{B_{2R_{1}} - B_{R_{1}}}\!\!\!\!\! 
                                    |\nabla u_{k}|^{p}| \xi_{k}|^{p} 
                                dxdy
                          \right)^{1/p} \!\! + 
                          \left(\int_{B_{2R_{1}} - B_{R_{1}}}\!\!\!\!\! 
                                |u_{k}|^{p}|\nabla \xi_{k}|^{p} 
                                dxdy
                          \right)^{1/p} \\
               & \leq & 
                          \left(\int_{B_{R_{k}} - B_{R_{0}}}\!\!\!\!\! 
                                     |\nabla u_{k}|^{p} 
                                     dxdy
                          \right)^{1/p} \!\! + 
                          \left(\int_{B_{2R_{1}} - B_{R_{0}}} \!\!\!\!\!
                                |\nabla \xi_{k}|^{p}| u_{k}|^{p} 
                                    dxdy
                          \right)^{1/p} \\
               &  <   & 
                           \varepsilon^{1/p} + 
                          \left(\int_{{\bf R}^{n+1}} \!\!\!\!\!
                                    |\nabla \xi_{k}|^{n+1} 
                                    dxdy
                          \right)^{\frac{1}{n+1}}
                          \left(\int_{B_{R_{k}} - B_{R_{0}}} \!\!\!\!\!
                                     |u_{k}|^{\frac{n+1}{n}q} 
                                     dxdy
                          \right)
                                ^{\frac{n}{(n+1)q}} \\ 
               &  <   & \varepsilon^{1/p} + C\varepsilon^{\frac{n}{(n+1)q}}.
\end{eqnarray*}
(All balls in the above are centered at $(w_{k}, \tilde{w}_{k})$.)
Similarly, we can show that 
$\;{\bf M}_{3}^{1/p} < \varepsilon + C \varepsilon^{\frac{n}{(n+1)q}}$. 
Denote 
$u_{1k} \equiv u_{k} \xi_{k}$, $u_{2k} \equiv u_{k}\eta_{k}$. 
By combining these estimates, we finally have
$$|{\bf M}|\! =\! 
\left| 
   \int_{{\bf R}^{n+1}} \!\!\!\! |\nabla u_{k}|^{p} dxdy \!-\!\! 
   \int_{{\bf R}^{n+1}} \!\!\!\! |\nabla u_{1k}|^{p} dxdy \!-\!\! 
   \int_{{\bf R}^{n+1}} \!\!\!\! |\nabla u_{2k}|^{p} dxdy 
\right|\!\! < \varepsilon + C \varepsilon^{\frac{np}{(n+1)q}}.$$ 
In other words, 
\[ 
    {\bf I} = \lim_{k \to \infty} \!
    \int_{{\bf R}^{n+1}}\!\!\! 
      |\nabla u_{k}|^{p} dxdy 
            =
              \lim_{k \to \infty} \!
                 \int_{{\bf R}^{n+1}}\!\!\! 
                 |\nabla u_{1k}|^{p} dxdy + 
              \lim_{k \to \infty} \!               
                 \int_{{\bf R}^{n+1}}\!\!\! 
                 |\nabla u_{2k}|^{p} dxdy. 
\]
It follows from the assumptions at the beginning that 
\begin{eqnarray}
&      & \left|\; \int_{{\bf R}^{n}} |u_{2k}|^{q} dx - 
           \left(\int_{{\bf R}^{n}} |u_{k}|^{q} dx - 
                 \int_{{\bf R}^{n}} |u_{1k}|^{q} dx 
           \right) \;
         \right| \nonumber \\
& \leq &  \int_{B_{R_{k}((w_{k}, \tilde{w}_{k}))} -
                B_{R_{1}((w_{k}, \tilde{w}_{k}))}} 
        |u_{k}(x,0)|^{q} \otimes \delta_{0}(y) dxdy \nonumber \\  
& \leq &   \int_{R_{0} \leq |(x,y) - (w_{k}, \tilde{w}_{k})| 
             \leq R_{k}} |u_{k}(x,y)|^{q} \otimes \delta_{0}(y) dxdy
          < \varepsilon. \label{**0}
\end{eqnarray}
Let 
$\alpha_{k} \equiv \int_{{\bf R}^{n}} |u_{1k}(x,0)|^{q} dx$, and  
$\beta_{k} \equiv \int_{{\bf R}^{n}} |u_{2k}(x,0)|^{q} dx $. 
By taking a subsequence, if necessary, we may assume that 
$\alpha_{k} \to \alpha$, and $\beta_{k} \to \beta$. We can see that
\vspace{-.1in}
\[
0 \leq  \alpha ,\; \beta \leq 1\;\; \mbox{and} \; \;
|\;\beta - (1 - \alpha)\;| < \varepsilon.
\] 
Use the estimates for ${\bf M}$ to observe that 
\[
\left|\; \int_{{\bf R}^{n+1}}\!\!\!\!\!
                 |\nabla u_{1k}(x,y)|^{p} 
           \! + \! |u_{1k}(x,y)|^{\frac{(n+1)q}{n}} 
           \! + \! |u_{1k}(x,0)|^{q} \otimes \delta_{0}(y) 
                      dxdy - \lambda^{*} \;\right|  <  
                                 \varepsilon, 
\] 
\[                                  
\left|\;\int_{{\bf R}^{n+1}} \!\!\!\!\!
                 |\nabla u_{2k}(x,y)|^{p} \!
           \! + \! |u_{2k}(x,y)|^{\frac{(n+1)q}{n}} 
           \! + \! |u_{2k}(x,0)|^{q} \!
                 \otimes \! \delta_{0}(y) 
                 dxdy - (L - \lambda^{*})  \right| <  \varepsilon. 
\]
We can also see that 
$\int_{{\bf R}^{n+1}} |\nabla u_{ik}(x,y)|^{p} dxdy 
\geq \gamma > 0$ for $i = 1,2$, and $\gamma$ a positive constant 
using the Sobolev embedding theorem 
and the Sobolev trace inequalities together with the estimates above.  
Now we look at all the possible values for 
$\alpha$ and $\beta$. They are: 
\begin{eqnarray*}
 (a) : \alpha_{k} \rightarrow 0 (\beta_{k} \rightarrow 1),  
& & (b) : \alpha \not= 0 (\beta \not= 1), \\ 
 (c) : \alpha_{k} \rightarrow 1 (\beta_{k} \rightarrow 0),  
& & (d) : \beta \not= 0 (\alpha \not= 1).  
\end{eqnarray*}
By exchanging the roles of $\alpha_{k}$ and $\alpha$  
with $\beta_{k}$ and $\beta$, 
the cases $(c)$ and $(d)$ reduce to the cases 
$(a)$ and $(b)$. 
In the case $(a)$, it follows from the estimates 
for ${\bf M}$ that 
${\bf I} \geq \gamma + {\bf I} - 
\varepsilon$ for all small $\varepsilon$, 
which leads to a contradiction that  
${\bf I} \geq \gamma + {\bf I} > {\bf I}$. 
For the case $(b)$, we define ${\bf I}_{\alpha}$ as 
\[
{\bf I}_{\alpha}\!\! \equiv  
\inf \!\: \left\{
         J(u) \equiv \! \int_{{\bf R}^{n+1}} 
 \!\!\! |\nabla u(x,y)|^{p} dxdy  :  
     \!\!  \int_{{\bf R}^{n}} \!\!\! |u(x,0)|^{q} dx 
              = \alpha, 
    u \in W^{1,p}({\bf R}^{n+1}) \right \}.              
\]
It easily follows from the definition that
${\bf I} = {\bf I}_{1}$ and 
${\bf I}_{\alpha} =  \alpha^{p/q}{\bf I}$. 
It can be also shown that 
\[
{\bf I} < {\bf I}_{\alpha} + 
                {\bf I}_{1 - \alpha} \;\; 
\mbox{for}\;\; 0 < \alpha < 1.
\]
This is called {\it Strict Subadditivity}.
Now, in the case $(b)$, we have 
${\bf I} \geq {\bf I}_{\alpha} + 
{\bf I}_{1 - \alpha} - \varepsilon $ 
for all small $\varepsilon > 0$ 
which violates the strict subadditivity.
This completes the proof of this lemma. 
\hfill$\Box$\par

\vspace{.13in}

Since we have shown that vanishing and 
dichotomy can not occur, 
we now conclude by Lemma \ref{con-1} that we have the compactness 
as follows:
there exists a sequence $(w_{k}, \tilde{w}_{k}) \in 
{\bf R}^{n} \times {\bf R}$ 
so that for any 
$\varepsilon > 0$, there is $R \in (0, \infty)$ 
such that 
\begin{eqnarray}
    \int_{[B_{R}((w_{k}, \tilde{w}_{k}))]^{C}}
         |\nabla u_{k}(x,y)|^{p} dxdy 
   +  \int_{[B_{R}((w_{k}, \tilde{w}_{k}))]^{C}} 
        |u_{k}(x,y)|^{\frac{(n+1)q}{n}} dxdy \nonumber \\
   +  
    \int_{[B_{R}((w_{k}, \tilde{w}_{k}))]^{C} \cap {\bf R}^{n} 
    \times \{ 0 \}}
    |u_{k} (x,0)|^{q} dx \;\;
    < \; \varepsilon.\;\; \label{**}
\end{eqnarray}
\smallskip

\begin{remark}
We may choose $\tilde{w}_{k} = 0$.
\end{remark}

\noindent {\bf Proof :}
If $\varepsilon < 1$, then $| \tilde{w}_{k}| \leq R$. 
[Otherwise, we would have 
$B_{R}((w_{k}, \tilde{w}_{k}))  \subset {\bf R}^{n} 
             \times ({\bf R} - \{ 0 \})$, 
implying 
\[
\int_{{\bf R}^{n}} \!\!\!
|u_{k}(x,0)|^{q}dx \;
\leq \;
\int_{[B_{R}((w_{k}, \tilde{w}_{k}))]^{C} \cap {\bf R}^{n} 
\times \{ 0 \} } \!\!\! P_{k}(x,y)dxdy \; < \; \varepsilon < 1,
\]
which violates the assumption that
$\int_{{\bf R}^{n}} |u_{k}(x,0)|^{q}dx = 1$.]
Take $(w_{k}, 0) \in {\bf R}^{n} 
\times \{ 0 \}$ and replace 
$R$ by $2R$. Then we have the compactness we had before
\[
\hspace{.5in}
\int_{[B_{2R}(w_{k}, 0)]^{C}} \!\!\! P_{k}(x,y)\:dxdy \;\leq \; 
\int_{[B_{R}(w_{k}, \tilde{w}_{k})]^{C}} \!\!\! P_{k}(x,y)\:dxdy \; 
< \varepsilon. \hspace{.5in} \mbox{\hfill$\Box$\par}
\]

\vspace{.13in}

We denote by $(u_{k})$ the new minimizing 
sequence $(\tilde{u}_{k})$ defined by 
$\tilde{u}_{k}(x, y) \equiv  u_{k}(x + w_{k}, y)$    
for all $(x, y) \in {\bf R}^{n} \times {\bf R}$. 
We may assume that 
\begin{eqnarray*}
u_{k} &\!\! \rightarrow u \;\;\;\; 
          \mbox{a.e. in}\;\; 
           {\bf R}^{n+1}, \;\; 
      &  u_{k}  \rightarrow u \;\; 
          \mbox{a.e. in}\;\; 
          {\bf R}^{n} \times 
                            \{ 0 \} \\
u_{k} & \rightharpoonup u \;\; 
          \mbox{in}\;\; W^{1,p}
           ({\bf R}^{n+1}),\;\;
      &  u_{k}  \rightharpoonup  u \;\; 
          \mbox{in}\;\; L^{q}
              ({\bf R}^{n}\times 
                           \{ 0 \}).
\end{eqnarray*}
\begin{lemma}[\hspace{-.015in}Concentration \hspace{-.03in}Compactness 
\hspace{-.03in}I\hspace{-.01in}I\hspace{-.015in}I\hspace{-.01in}] \label{con-2}
Let $(u_{k})$ be a bounded sequence in 
$W^{1,p}({\bf R}^{n+1})$ such that 
$(|\nabla u_{k}(x,y)|^{p})$ is tight. We may assume 
$u_{k} \rightarrow u$ a.e. in ${\bf R}^{n+1}$ 
and $(|\nabla u_{k}(x,y)|^{p})$ and 
$(|u_{k}(x,0)|^{q} \otimes \delta_{0}(y))$ converge 
weakly to some bounded nonnegative measures  $\mu, \nu$ on 
${\bf R}^{n+1}$
and supp $(\nu) \subset {\bf R}^{n} \times \{ 0 \}$. Then \\
$(i)$ There exist some at most countable set ${\cal L}$ and 
two families $(x_{l})_{l \in {\cal L}}$ of distinct points in 
${\bf R}^{n}, \; (\nu_{l})_{l \in {\cal L}}$ in 
$(0, \infty)$ such that 
\begin{eqnarray*}
    \nu &   =   & |u(x,0)|^{q} \otimes \delta_{0}(y) 
            + \sum_{l \in {\cal L}} \nu_{l} \delta_{(x_{l},0)} \\
    \mu & \geq  & |\nabla u(x,y)|^{p} 
            + \sum_{l \in {\cal L}} {\bf I}\; 
              \nu_{l} \delta_{(x_{l},0)}.
\end{eqnarray*}
$(ii)$ If $u = 0$ and $\mu({\bf R}^{n+1}) \leq 
{\bf I}\; \nu({\bf R}^{n+1})^{p/q}$,   
then ${\cal L}$ is a singleton and  
$\nu = c_{0}\; \delta_{(x_{0},0)}$, and
$\mu = {\bf I}\;c_{0}^{p/q} 
                      \;\delta_{(x_{0},0)}$ 
for some $c_{0} > 0$, and for some $x_{0} \in {\bf R}^{n}$.
\end{lemma}

\noindent {\bf Proof :} 
We first look at the case $u \equiv 0$. By the Sobolev trace 
inequality, we have for $\varphi \in C_{0}^{\infty}({\bf R}^{n+1})$ 
\begin{eqnarray}
\left( 
 \!\! \int_{{\bf R}^{n}} \!\!\! |\varphi(x,0)u_{k}(x,0)|^{q} dx
\right)^{1/q} \!\!\!\!
\leq {\bf I}^{-1/p} \!\!
\left(
  \int_{{\bf R}^{n+1}} \!\!\! 
      |\nabla [\varphi(x,y)u_{k}(x,y)]|^{p} dxdy 
\right)^{1/p}.  \label{tracein}
\end{eqnarray}
The left-hand side of (\ref{tracein}) 
\[
\left( 
 \int_{{\bf R}^{n}}\!\! |\varphi(x,0)u_{k}(x,0)|^{q} dx
\right)^{1/q} \!\!= 
\left( 
 \int_{{\bf R}^{n+1}}\!\! |\varphi(x,y)u_{k}(x,0)|^{q} 
                                \otimes \delta_{0}(y) dxdy
\right)^{1/q}
\]
converges to 
$\left( 
  \int_{{\bf R}^{n+1}} |\varphi|^{q} d\nu
\right)^{1/q}$ 
as $k \to \infty$. By the fact that  $\varphi$ has compact 
support, $(u_{k})$ converges to $0$ a.e., and 
\begin{eqnarray*} 
\left| 
\| \nabla (\varphi u_{k}) \|_{L^{p}({\bf R}^{n+1})} - 
\| \varphi (\nabla u_{k}) \|_{L^{p}({\bf R}^{n+1})}
\right|^{p} 
\!\!& \leq &\!\!
\int_{{\bf R}^{n+1}} \!\!|\nabla \varphi(x,y)|^{p}|u_{k}(x,y)|^{p} dxdy
\end{eqnarray*}
converges to $0$ as $k \to \infty$, we have 
that the right-hand side of (\ref{tracein}) converges to 
$\left(
  \int_{{\bf R}^{n+1}} |\varphi(x,y)|^{p} d \mu
 \right)^{1/p}$.
Taking $k \to \infty$ in (\ref{tracein}) yields for 
$\varphi \in C_{0}^{\infty}({\bf R}^{n+1})$,  
\[
\left(
  \int_{{\bf R}^{n+1}} |\varphi(x,y)|^{q} d \nu
\right)^{1/q}
\leq {\bf I}^{-1/p} 
\left(
  \int_{{\bf R}^{n+1}} |\varphi(x,y)|^{p} d \mu 
\right)^{1/p}.
\] 
By applying Lemma \ref{con-3} to two measures 
$\mu$ and $\nu$ on ${\bf R}^{n+1}$, 
we obtain the results of Lemma \ref{con-2}. 
Now consider the general case that the weak limit $u$ 
is not necessarily $0$. Let $v_{k} = u_{k} - u$. 
By applying to $(v_{k})$ what we have proved for $(u_{k})$ above 
and using Br\'ezis-Lieb lemma saying that for 
$\varphi \in C_{0}^{\infty}({\bf R}^{n})$ 
\[
\int_{{\bf R}^{n}} |\varphi|^{q}|u_{k}|^{q} dx -
\int_{{\bf R}^{n}} |\varphi|^{q}|v_{k}|^{q} dx  
\;\; \mbox{converges to}\;\;
\int_{{\bf R}^{n}} |\varphi|^{q}|u|^{q} dx, 
\]
we have the representation for 
$\nu = |u(x,0)|^{q} \otimes \delta_{0}(y) + 
 \sum_{l \in {\cal L}} \nu_{l} \delta_{(x_{l},0)}$ for some 
countable set ${\cal L}$. 
We have that for $\varphi \in C_{0}^{\infty}({\bf R}^{n+1})$, 
\vspace{-.1in}
\begin{eqnarray}
\hspace{-.13in}&      &\hspace{-.08in}
{\bf I}^{1/p}  
\left(
  \int_{{\bf R}^{n+1}} 
     |\varphi(x,y)|^{q}|u_{k}(x,y)|^{q} \otimes \delta_{0}(y) dxdy
\right)^{1/q} \nonumber  \\  
\hspace{-.13in}& \leq &\hspace{-.08in}\!\!
\left(\!
  \int_{{\bf R}^{n+1}} \!\!\!\!\!\!
     |\varphi(x,y)|^{p}|\nabla u_{k}(x,y)|^{p} dxdy \!
\right)^{1/p}  
\!\!\!\!\! + \!
\left(\!
  \int_{{\bf R}^{n+1}} \!\!\!\!\!\!
     |\nabla \varphi(x,y)|^{p}|u_{k}(x,y)|^{p} dxdy \!
\right)^{1/p}  \label{trace2}
\end{eqnarray}
and 
\vspace{-.12in}
\[
\int_{{\bf R}^{n+1}} |\nabla \varphi|^{p}|u_{k}|^{p} dxdy 
\;\;\;\mbox{converges to} \;\;\; 
\int_{{\bf R}^{n+1}} |\nabla \varphi|^{p}|u|^{p} dxdy,
\]
since $|\nabla \varphi|$ has compact support. 
Passing to the limit in (\ref{trace2}), we have 
\begin{eqnarray}
{\bf I}^{1/p} \!\! 
\left(
   \int_{{\bf R}^{n+1}} \!\!\! |\varphi|^{q} d \nu
\right)^{1/q} \!\!\!\!
\leq \!\!
\left(
   \int_{{\bf R}^{n+1}} \!\!\! |\varphi|^{p} d \mu
\right)^{1/p} \!\!\!\! + \!
\left(
   \int_{{\bf R}^{n+1}} \!\!\!
       |\nabla \varphi|^{p}|u|^{p} dxdy
\right)^{1/p}. \label{trace3}
\end{eqnarray}
Take $\varphi \in C_{0}^{\infty}({\bf R}^{n+1})$ 
satisfying $0 \leq \varphi \leq 1$, $\varphi(0) = 1$, and 
$\mbox{supp}\; \varphi = B_{1}(0)$. 
Apply (\ref{trace3}) to 
$\varphi(\frac{x-x_{l}}{\varepsilon},\frac{y}{\varepsilon})$,  
for $l \in {\cal L}$ and $\varepsilon$ positive 
and small enough, to have  
\begin{eqnarray*}
&      &
{\bf I}^{1/p} 
\nu \left( B_{\varepsilon}(x_{l},0) \right)^{1/q} \\ 
& \leq & \mu \left( B_{\varepsilon}(x_{l},0)
             \right)^{1/p} 
       + \left( 
           \int_{B_{\varepsilon}(x_{l},0)} 
              |\nabla [ \varphi
                              (\frac{x-x_{l}}{\varepsilon},
                               \frac{y}{\varepsilon})]
                      |^{p}|u(x,y)|^{p} dxdy
         \right)^{1/p}. 
\end{eqnarray*} 
By the Sobolev embedding theorem, we have 
\begin{eqnarray*}
&      &
\left( 
    \int_{B_{\varepsilon}(x_{l},0)}
          |\nabla [ \varphi
                              (\frac{x-x_{l}}{\varepsilon},
                               \frac{y}{\varepsilon})]
                      |^{p}|u(x,y)|^{p} dxdy
\right)^{1/p} \\
& \leq &
\left(
     \int_{B_{\varepsilon}(x_{l},0)}\!\!\!\! 
         |u(x,y)|^{\frac{n+1}{n}q} dxdy
\right)^{\frac{n}{(n+1)q}} \!\!\!
\left(
     \int_{B_{\varepsilon}(x_{l},0)}\!
         \left| \nabla [\varphi(\frac{x-x_{l}}{\varepsilon},
                                    \frac{y}{\varepsilon})]
         \right|^{n+1} \!\!\!\!\!\!\! dxdy 
\right)^{\frac{1}{n+1}} \\
& \leq & 
\left(
    \int_{B_{\varepsilon}(x_{l},0)}\!\!\!\! 
           |u(x,y)|^{\frac{n+1}{n}q} dxdy
\right)^{\frac{n}{(n+1)q}}\!\! 
\left(
    \int_{{\bf R}^{n+1}}\!\!\!\! 
           |\nabla \varphi(x,y)|^{n+1} dxdy 
\right)^{1/(n+1)} \\
& \leq & 
D \left( 
    \int_{B_{\varepsilon}(x_{l},0)} 
           |u(x,y)|^{\frac{n+1}{n}q} dxdy
\right)^{\frac{n}{(n+1)q}} \label{trace4}
\end{eqnarray*}
where $D$ is a positive constant.
Taking $\varepsilon \to \infty$ yields  
\begin{eqnarray*}
{\bf I}^{1/p} \nu(\{ (x_{l},0) \})^{1/q} & \leq & 
                    \mu(\{ (x_{l},0) \})^{1/p}, \\ 
\mbox{then}\;\;\;\; {\bf I}^{1/p} \nu_{l}^{1/q} & \leq & 
                     \mu(\{ (x_{l},0) \})^{1/p},  \\ 
\mbox{and so},\;\;\;\;\;\;\;\;\;    \mu & \geq & 
 {\bf I}\; \nu_{l}^{p/q} \delta_{(x_{l},0)} \;\;
\mbox{for}\;\;\; l \in {\cal L}. 
\end{eqnarray*}
Thus, 
$\mu \geq \sum_{l \in {\cal L}} {\bf I} 
                  \nu_{l}^{p/q} \delta_{(x_{l},0)}$. 
Let $\sum_{l \in {\cal L}} {\bf I} 
                  \nu_{l}^{p/q} \delta_{(x_{l},0)} = \mu_{1}$. 
By the fact that two measures $\mu_{1}$ and $|\nabla u|^{p}$ are 
orthogonal, and $\mu \geq |\nabla u|^{p}$ by the weak convergence, 
we conclude that  
$\mu \geq |\nabla u|^{p} + 
         \sum_{l \in {\cal L}} {\bf I} 
                  \nu_{l}^{p/q} \delta_{(x_{l},0)}$ 
to complete the proof.
\hfill$\Box$\par

\begin{lemma} 
$u \not \equiv 0.$ \label{notequal} 
\end{lemma}

\noindent {\bf Proof :}
Suppose $u \equiv 0$. 
Then $(u_{k})$ converges weakly to $0$ in 
$W^{1,p}({\bf R}^{n+1})$. 
We know that $(|\nabla u_{k}(x,y)|^{p})$ converges weakly to $\mu$ 
tightly in the space of measures and 
$(|u_{k}(x,0)|^{q} \otimes \delta_{0}(y))$ 
converges weakly to $\nu \;\;(\;\mbox{supp}(\nu) 
\subset {\bf R}^{n} \times \{ 0 \})$ from $(\ref{**})$.  
We can see that 
\begin{eqnarray*}
   \int_{{\bf R}^{n+1}} d\mu 
   & = & \lim_{k \to \infty} 
   \int_{{\bf R}^{n+1}} 
     |\nabla u_{k}(x,y)|^{p} dxdy 
    = {\bf I}, \\ 
   \int_{{\bf R}^{n+1}} d\nu \; 
   & = & \lim_{k \to \infty} 
   \int_{{\bf R}^{n+1}} |u_{k}(x,0)|^{q} 
      \otimes \delta_{0}(y) dxdy \; 
    = \; 1.
\end{eqnarray*}  
In other words, 
\vspace{-.1in}
\[ 
   \mu({\bf R}^{n+1}) \; = \; {\bf I} 
   \;  = \; {\bf I} \;\nu({\bf R}^{n+1}).  
\]
By Lemma \ref{con-3}, there exists 
$x_{0} \in {\bf R}^{n}$, and so that $\nu = \delta_{(x_{0},0)}$ and 
$\mu = {\bf I} \delta_{(x_{0},0)}$. 
Then it gives 
a contradiction saying
\[
\frac{1}{2} = Q_{k}(1) \geq 
 \int_{B_1(x_{0},0)} \! |u_{k}(x,0)|^{q} 
 \otimes \delta_{0}(y) dxdy 
 \rightarrow \! \int_{B_1(x_{0},0)} \! d\nu = 
 \nu(B_1(x_{0},0)) = 1,
\] 
so we complete the proof.   
\hfill$\Box$\par

Let $\int_{{\bf R}^{n}} |u(x,0)|^{q}dx =   
\int_{{\bf R}^{n+1}} |u(x,0)|^{q} \otimes 
\delta_{0}(y) dxdy = \alpha$. From Lemma \ref{notequal}, we have 
$0 < \alpha \leq 1$. Now it is sufficient to show that 
$\alpha = 1$ in order to prove Theorem \ref{2.2.1}. 
So suppose $\alpha \not = 1$. By Lemma \ref{con-2}, 
there exist a set ${\cal L}$ at most countable, 
$ (x_{l})_{l \in {\cal L}} \subset {\bf R}^{n}$ and 
$(\nu_{l})_{l \in {\cal L}} \in (0, \infty)$ such that 
\begin{eqnarray*}
\nu &  =     & |u(x,0)|^{q} \otimes \delta_{0}(y) + 
                \sum_{l \in {\cal L}} \nu_{l}\delta_{({x_{l}},0)} 
\;\; (\mbox{so,}\;\; 1 = \alpha + \sum_{l \in {\cal L}} \nu_{l} ), \\
\mu &  \geq  & |\nabla u(x,y)|^{p} + 
              \sum_{l \in {\cal L}}{\bf I}
              \nu_{l}^{p/q}\delta_{({x_{l}},0)}. 
\end{eqnarray*}
This leads us to a following contradiction:
\begin{eqnarray*}
{\bf I}_{\alpha} 
          & \leq & \int_{{\bf R}^{n+1}}
                                  |\nabla u(x,y)|^{p}dxdy\\
          & \leq & \int_{{\bf R}^{n+1}}
                                 d\mu - 
                                \int_{{\bf R}^{n+1}} 
                                  \sum_{l \in {\cal L}}{\bf I}
                                  \nu_{l}^{p/q}\delta_{({x_{l}},0)}
                                 dxdy\\
          &  =  & {\bf I} - \sum_{l \in {\cal L}}
                                  {\bf I}\nu_{l}^{p/q} 
             =   {\bf I}(1 - \sum_{l \in {\cal L}}
                                  \nu_{l}^{p/q})\\
          &  <  & {\bf I}(1 - \sum_{l \in {\cal L}}
                                  \nu_{l})^{p/q}
             =    {\bf I}\;\alpha^{p/q} = 
                               {\bf I}_{\alpha}. 
\end{eqnarray*}
The last inequality holds since 
$\sum_{l \in {\cal L}} \nu_{l} = 1 - \alpha \neq 0$. 
So we conclude that $\alpha = 1$, and this proves that 
there exists an extremal function for the trace inequality. 
This completes the proof of Theorem \ref{2.2.1}. 
\hfill$\Box$\par

\section{Conjectured extremal function}                            

\indent

Any extremal function for the Sobolev trace inequality satisfies 
the following Euler-Lagrange equation: for a positive constant $C$,
\vspace{.1in}
\begin{eqnarray}
    \left. \begin{array}{ll}
    \mbox{div}(|\nabla u|^{p-2}\nabla u) = 0 & \hspace{1in} \mbox{on}\;\; 
                                  {\bf R}^{n+1}_{+}        \\
                            \\ 
      |\nabla u|^{p-2} \frac{\partial u}{\partial y} +
               C |u|^{q-2}u      = 0 & \hspace{1in} \mbox{on}\;\; 
                         \partial {\bf R}^{n+1}_{+} 
        \end{array} \right\}.  \label{*4} 
\end{eqnarray}

\noindent
It can be easily proved that there is no radial function in all 
the variables in $(x,y) \in {\bf R}^{n+1}_{+}$ satisfying 
the equation (\ref{*4}) due to the boundary condition. 
As a way to identify a function which satisfies (\ref{*4}), we
will look at a restricted class of functions. In particular, 
we will restrict our attention to the functions of conformal 
factor, $[(1+y)^{2}+|x|^{2}]$, where $(x,y) \in {\bf R}^{n+1}_{+}$. 
This means we assume an extra 
symmetry for possible extremal functions for the Sobolev trace 
inequality. This choice of symmetry is not surprising 
if we look at the extremal function for the special case of 
the Sobolev trace inequality with $p = 2$. This choice also 
specifies the function on the boundary as a function of 
$[1+|x|^{2}]$. This condition is not at all strict since 
it suffices to consider radial decreasing functions on 
${\bf R}^{n}$ for extremal functions by using a rearrangement 
technique. J. Escobar conjectured the extremal function for 
the Sobolev trace inequality in \cite{Escobar} as 
$[(1+y)^{2}+|x|^{2}]^{- \frac{n+1-p}{2(p-1)}}$. 
The following remark will make it clear 
that it is the only possible choice of function for 
the extremal function.

\begin{remark}[\hspace{-.011in}the \hspace{-.02in}conjectured 
\hspace{-.02in}extremal \hspace{-.02in}for \hspace{-.02in}the 
\hspace{-.02in}Sobolev \hspace{-.02in}trace 
\hspace{-.02in}inequality\hspace{-.011in}] 
Suppose that $f$ is an extremal function for the Sobolev 
trace inequality and that $f$ is a function of $[(1+y)^{2}+|x|^{2}]$. 
We may also assume that $f$ is decreasing in $|x|$, and in $y$. 
Then $f(x,y)$ is exactly the same function that was conjectured. 
\end{remark}
\noindent {\bf Proof :}
Let $f(x,y) \equiv  \Phi (v(x,y))$, 
where $\Phi$ is a function of one variable and 
$v(x,y) \equiv  (1+y)^{2}+|x|^{2}$. 
This gives the following equations:
\begin{eqnarray*}
\frac{\partial f}{\partial x_{j}}(x,y) &\!\!\!\!\!\!\!\!\!\!\!\!\!
  \!\!\!\!\!\!\!\!\!\!\!\!\!\!\!\!\!= 
2x_{j} \Phi'(v)
      &\:\:\frac{\partial f}{\partial y}(x,y) =  2(1+y) \Phi'(v) \\
\frac{\partial^{2} f}{\partial x_{j}^{2}}(x,y) & \!\!=   
                       2[2x_{j}^{2} \Phi''(v) + \Phi'(v)]
      &\frac{\partial^{2} f}{\partial y^{2}}(x,y)      =  
                       2[2(1+y)^{2} \Phi''(v) + \Phi'(v)],  
\end{eqnarray*}
where $'$ denotes the derivative with respect to $v$.
These equations and the fact that $f$ satisfies (\ref{*4}) 
since $f$ is an extremal function yield 
\[
\mbox{div}(|\nabla f|^{p-2} \nabla f) = 
    2^{p-1}|\Phi'|^{p-2} v^{\frac{p}{2} - 1}
    [2(p-1) \Phi''(v)v + (n+p-1) \Phi'(v)] = 0.
\]
Since $f$ is not a constant function,
we have the equation that $\Phi$ must satisfy:
\[ 
[2(p-1) \Phi''(v)v + (n+p-1) \Phi'(v)] = 0.
\]
From this, we have 
\[
\left[ \ln |\Phi'(v)| \right]'
    = \frac{\Phi''(v)}{\Phi'(v)} 
    = - \frac{n+p-1}{2(p-1)} \frac{1}{v}.
\]
Hence we obtain 
$\Phi(v) = c_{0} v^{- \frac{n+1-p}{2(p-1)}} = 
c_{0} [(1+y)^{2}+|x|^{2}]^{- \frac{n+1-p}{2(p-1)}}$, for some 
constant $c_{0}$, which can be determined uniquely by 
the condition that $\| f \|_{L^{q}({\bf R}^{n})} = 1$. 
This function is the very function that Escobar conjectured.
\hfill$\Box$\par 

\vspace{.2in}

The following proposition characterizes this function as  
the minimizer of the Sobolev trace inequality functional 
when restricted 
to the class of functions of conformal factor. For this we define 
\[
{\cal J}(\omega) \equiv  \int_{{\bf R}^{n+1}_{+}}    
                       | \nabla \omega(x,y) |^{p}dxdy,  
\]  
where $\omega$ belongs to the admissible set 
\begin{eqnarray*}
{\cal A} \equiv \{ \omega \in W^{1,p}({\bf R}^{n+1}_{+}) : 
&&\hspace{-.2in} \omega \mbox{ is a function of } 
[(1+y)^{2}+|x|^{2}], \mbox{ and } \\
&&\hspace{-.2in} \omega(x,0) = 
c_{0}(1 + |x|^{2})^{{- \frac{n+1-p}{2(p-1)}}}  
                 \}.
\end{eqnarray*}
\begin{proposition}
Let $f$ be the conjectured extremal function for 
the Sobolev trace inequality. Then 
\vspace{-.1in}
\[
{\cal J}(f) = \min_{\omega \in {\cal A}}{\cal J}(\omega),
\]
in other words, the infimum of $\cal J(\cdot)$ on $\cal A$
is attained at $f$.
\end{proposition}
\noindent {\bf Proof :}
Take any $\omega \in {\cal A}$ and consider $f- \omega$. 
Since $f$ satisfies the equation (\ref{*4}), 
we have  
\[
\mbox{div} \left( |\nabla f|^{p - 2}\nabla f \right) 
\left(f- \omega \right) = 0 \;\; \mbox{on} 
                 \;\; {\bf R}^{n+1}_{+}.        
\]
An integration by parts yields 
\[
0  =  \int_{{\bf R}^{n+1}_{+}}
        | \nabla f |^{p-2}\left( |\nabla f|^{2}
                 -  \nabla f \cdot \nabla\omega \right) dxdy 
\]
and there is no boundary term since  
$f - \omega = 0$ on $\partial {\bf R}^{n+1}_{+}$
by the fact that both $f$ and $\omega$ belong to ${\cal A}$. 
Now Young's inequality gives 
\begin{eqnarray*}
{\cal J}(f) & = & \int_{{\bf R}^{n+1}_{+}}    
                       | \nabla f |^{p}dxdy 
             =  \int_{{\bf R}^{n+1}_{+}}    
                       | \nabla f |^{p-2}
                \left( \nabla f \cdot \nabla \omega \right) dxdy \\
            & \leq & (1 - \frac{1}{p})\int_{{\bf R}^{n+1}_{+}} 
            |\nabla f|^{p}dxdy + 
   \frac{1}{p} \int_{{\bf R}^{n+1}_{+}}|\nabla\omega|^{p}dxdy.              
\end{eqnarray*} 
We obtain 
\vspace{-.1in}
\[
\hspace{1.302in}
{\cal J}(f) \leq {\cal J}(\omega) \hspace{.7in}
(\omega \in {\cal A}).
\hspace{1.302in}\mbox{\hfill$\Box$\par}
\]

\section{Sobolev trace inequality with $\mathbf{p=1}$}   
\indent

In this section, we will treat the Sobolev trace inequality 
for the case when $p=1$ (thus $q=1$) separately. 
The existence of the extremal function for the Sobolev trace 
inequality for 
the case when $p = 1$ ($q =1$) is not guaranteed by the argument 
used for $p$ with $1 < p < n+1$. This is one of the limit cases 
of the inequality and is closely related to the isoperimetric 
inequality. 

The Sobolev trace inequality for $p=1$ is given by  
\[
\int_{{\bf R}^{n}}| u(x,0) | dx \leq C 
\int_{{\bf R}^{n+1}_{+}}| \nabla u(x,y) |dxdy
\]
for a positive constant $C$. 
To find the best constant for this inequality, 
we look at the following quotient:
\vspace{-.1in}
\[
J(u) \equiv \frac{\left(\int_{{\bf R}^{n+1}_{+}}    
                       | \nabla u(x,y) |dxdy \right)}
            {\left(\int_{{\bf R}^{n}}   
                              | u(x,0) | dx \right)},
\] 
\noindent
where $u \in W^{1,1}({\bf R}^{n+1}_{+})$ and $u \not\equiv 0$. 
The best constant ${\bf I}$ is defined by 
\begin{eqnarray*}
{\bf I} \equiv  \inf  \{ J(u) : u \in W^{1,1}
    ({\bf R}^{n+1}_{+}), \;\;u \neq 0 \}. 
\end{eqnarray*}

Define ${\cal B}  \equiv 
     \left\{ g \in W^{1,1}({\bf R}^{n+1}_{+}) :  
                  g \geq 0 \;\;\mbox{on}\;\; {\bf R}^{n+1}_{+}, \; 
               \int_{{\bf R}^{n}} g(x,0) \: dx = 1 
     \right\}$. 
It is sufficient to consider functions in ${\cal B}$ to compute 
the best constant, since $J(\cdot)$ is dilation invariant and 
$J(u) = J(|u|)$.
Moreover, we will use a rearrangement technique to reduce further 
the functions to consider to 
a class of functions with a special property. Namely, we will take 
$\Phi^{*}_{S}$ to be the {\it Steiner rearrangement} of $\Phi$. Here 
$\Phi^{*}_{S}$ is symmetric radial decreasing in $x$, and is 
decreasing in $y$. 
Then we know that 
\begin{eqnarray}
\int_{{\bf R}^{n}} |\Phi(x,0)|dx  
   & = &  \int_{{\bf R}^{n}}|\Phi^{*}_{S}(x,0)|dx  \nonumber \\
   & = &  \int_{{\bf R}^{n}}\int_{0}^{\infty} 
             - \frac{\partial \Phi^{*}_{S}}{\partial y}(x,y) dydx  
        \nonumber\\
   & \leq &  \int_{{\bf R}^{n+1}_{+}}|\nabla \Phi^{*}_{S}(x,y)|dxdy
        \label{3.4}\\
   & \leq &  \int_{{\bf R}^{n+1}_{+}}|\nabla \Phi(x,y) |dxdy. 
             \label{3.5} 
\end{eqnarray} 
By the above observation, it suffices to consider 
functions in ${\cal B}$ having the following property $(P)$:

$(P)$: $g$ is symmetric radial decreasing in $x$ and 
decreasing in $y$.    

\noindent            
For any function $g$ having the property $(P)$, the  
inequality (\ref{3.5}) becomes equality.  
It is now clear that $\inf \{ J(g) \mid g\;\;
\mbox{has the property (P)}, \;\; g  \in {\cal B} \} \geq 1$. 

\begin{theorem}
\[
{\bf I} \equiv  \inf  \{ J(u) : u \in W^{1,1}
    ({\bf R}^{n+1}_{+}), \;\;u \neq 0 \} = 1.
\]
\end{theorem}

\noindent {\bf Proof :}
We will look at the inequalities above. 
The inequality $(\ref{3.5})$ becomes equality, since we choose $f$ 
with the property $(P)$. The question is when the inequality 
$(\ref{3.4})$ becomes equality. For that we require that $f$ satisfy 
\[
\left| \frac{\partial f}{\partial y}(x,y) \right| = 
|\nabla f(x,y)| \;\;\mbox{on} \;\;
{\bf R}^{n+1}_{+}.
\]
This means that 
\[
\frac{\partial f}{\partial x_{j}}(x,y) = 0\;\;
\mbox{on}\;\;{\bf R}^{n+1}_{+} \;\;\;\;\mbox{for} 
                  \;\;\; j = 1,2,3, \cdots ,n.
\] 
From this, we can see that $f$ should be a 
function of $y$ variable only. On the other hand, 
$f(x,0)$ is a function in 
$L^{p}({\bf R}^{n})$, so we need some restriction on 
the function. Any function of $y$ with appropriate decay 
multiplied by a characteristic 
function in the $x$ variable will be an extremal function. 
The problem is that such functions do not belong to 
$W^{1,1}({\bf R}^{n+1}_{+})$, which means that the extremal function 
does not exist. However, we can use an approximation argument to compute 
the best constant. Take a function 
$f(x,y) = \phi(y)\chi_{B}(x)$, where $\phi$ is a positive 
non-increasing function of $y$ variable 
and $B$ is the unit ball centered at the origin 
in ${\bf R}^{n}$. Then we have 
\[
\int_{{\bf R}^{n+1}_{+}}|\nabla f(x,y)|dxdy = 
\int_{{\bf R}^{n}} |f(x,0)| dx + 
\sigma_{n}\int_{0}^{\infty} \phi(y)dy 
\]
where $\sigma_{n}$ is the surface area of the unit ball in 
${\bf R}^{n}$. If we can make the second term 
in the right hand side go away, then we get the claim we made. 
Let $\phi_{\varepsilon}(y) = \exp(-\frac{\pi y^{2}}{\varepsilon})$. Then 
$\int_{0}^{\infty} \phi_{\varepsilon}(y)dy = \sqrt{\varepsilon}$, 
so that we can make it as small as we want.

\end{document}